\documentclass{article}

\usepackage{amssymb}
\usepackage[cp1251]{inputenc}
\usepackage[russian, english]{babel}
\usepackage{amsmath}
\usepackage{color}
\usepackage{graphics}
\usepackage{epsfig}
\usepackage[all,knot]{xy}
\xyoption{arc}

\def\mapr#1{\smash{\mathop{\buildrel{#1}\over\longrightarrow}}}

\newtheorem{theorem}{Theorem}

\def\proof{{\bf Proof.}}
\def\qed{\hfill\vrule width2mm height2mm depth2mm}

\def\A{{\bf A}}
\def\C{{\bf C}}
\def\F{{\bf F}}
\def\R{{\bf R}}
\def\T{{\bf T}}
\def\X{{\bf X}}
\def\cA{{\cal A\;}}
\def\cG{{\cal G}}
\def\cH{{\cal H}}

\def\Ad{\hbox{\rm Ad}}
\def\ad{{\hbox{\bf ad}}}
\def\Der{{\hbox{\bf Der}}}
\def\bydef{\stackrel{def}{=}}
\def\Int{{\hbox{\bf Int}\;}}
\def\Mor{\hbox{\rm Mor}}
\def\Obj{\hbox{\rm Obj}}
\def\Out{{\hbox{\bf Out}\;}}

\def\blf{$}
\def\elf{$}

\def\bdf{$$}
\def\edf{$$}

\def\beq#1{\begin{equation}\label{#1}}
\def\eeq{\end{equation}}

\def\mat{$}
\def\tam{$}

\newenvironment{dedication}
{
   \cleardoublepage
   \thispagestyle{empty}
   \vspace*{\stretch{1}}
   \hfill\begin{minipage}[t]{0.66\textwidth}
   \raggedright
} {
   \end{minipage}
   \vspace*{\stretch{3}}
   \clearpage
}

\newtheorem{proposition}{Proposition}

\title{
Derivations of Group Algebras\thanks{The text is presented in two languages - English and Russian}}
\author{Arutyunov A.~A.\thanks{Financially supported by Ministry of Education
and Science of the Russian Federation (Agreement no.~02.a03.21.0008 as of
24.06.2016).}, Mishchenko A.~S.\thanks{Partially supported by the RFBR grant
no.~14-01-00007.}, Shtern A.~I.\thanks{Partially supported by the RFBR grant
no.~14-01-00007.}}

\begin{document}
\begin{dedication}{Dedicated to the blessed memory of Yu.~P.~Solov'\"ev}\end{dedication}
\maketitle

\begin{abstract}
In the paper, a method of describing the outer derivations of the group algebra
of a finitely presentable group is given. The description of derivations is
given in terms of characters of the groupoid of the adjoint action of the
group.
\end{abstract}

\section{Introduction}

Consider an algebra \mat\cA\tam and some bimodule \mat E\tam over the algebra
\mat\cA\tam. Denote by \mat\Der(\cA,E)\tam the space of all derivations from
the algebra \mat\cA\tam to the bimodule \mat E\tam, i.e., the set of mappings
\bdf D\colon\cA\mapr{}E, \edf satisfying the condition \bdf D(ab)=D(a)b+aD(b),
\quad a,b\in\cA \edf (see~Losert(2008) \cite{Losert-2008}, Ghahramani(2000)
\cite{Ghahramani-2000}). Among the derivations \mat\Der(\cA,E),\tam we can
single out the so-called inner derivations \mat\Int(\cA,E)\subset\Der(\cA,E)
\tam that are defined by the adjoint representation, \bdf \ad_{x}(a)\bydef
xa-ax, \quad x\in E, a\in\cA. \edf

The derivation problem is formulated as follows: Are all derivations inner?
This problem was considered for group algebras \mat\cA=C[G]\tam of some
group~\mat G\tam rather than for all algebras. To be more precise, the group
algebra \mat \overline\cA=L^{1}(G)\tam and the bimodule \mat E=M(G)\tam is
considered, where \mat M(G)\tam stands for the algebra of all bounded measures
on~\mat G\tam with the multiplication operation defined by the convolution of
measures.

A question in the paper Dales(2000) \cite{Dales-2000}, (Question~5.6.B, p.~746)
is formulated as follows: Let \mat G\tam be a locally compact group. Does every
derivation from the algebra \mat \cA=L^{1}(G)\tam to the bimodule \mat
E=M(G)\tam be an inner derivation? An affirmative answer is supported by the
following consideration.

For the case in which the group \mat G\tam is a finitely generated discrete
free Abelian group, i.e., \mat G\approx\mathbb{Z}^{n}\mat, the algebra \mat
\overline\cA=L^{1}(G)\tam can also be identified with the Fourier algebra \mat
A(\mathbb{T}^{n})\tam of continuous functions on the \mat n\tam\!\!-dimensional
torus \mat\mathbb{T}^{n}\tam whose Fourier coefficients form an absolutely
convergent multiple series, \mat \cA =A(\mathbb{T}^{n})\subset
C(\mathbb{T}^{n}),\tam  (the Fourier algebra is smaller than the algebra of
continuous functions). There are no derivations on the algebra \mat
A(\mathbb{T}^{n})\tam because it has sufficiently many nonsmooth functions;
certainly, there are no inner derivations as well, because the algebra \mat
\overline\cA=L^{1}(G)\tam is commutative.

We are interested however in a dense subalgebra \mat
\cA=C[G]\subset\overline\cA\mat\ of the whole Banach algebra \mat
\overline\cA=L^{1}(G)\tam (rather than in the algebra \mat
\overline\cA=L^{1}(G)\tam itself); this subalgebra consists of the so-called
smooth elements of the algebra \mat \overline\cA=L^{1}(G),\tam in the
terminology of Connes (\cite{Connes-1994}, p.~247). For the group algebra \mat
\cA=C[G]\;\mat, one can also formulate a similar problem: Describe the algebra
of all outer derivations of the group algebra \mat \cA=C[G].\mat

\section{Group algebra $C[G]$.}

Consider the group algebra $\cA=C[G].$ Assume that \mat G\tam is a finitely
presentable discrete group.

An arbitrary element $u\in\cA$ is a finite linear combination
$u=\sum\limits_{g\in G}\lambda^{g}\cdot g. $ Consider an arbitrary linear
operator
$$
X\colon\cA\mapr{}\cA.
$$
A linear operator $X$ has the following matrix form: \beq{1}
X(u)=\sum\limits_{h\in G}\left(\sum\limits_{g\in
G}x_{g}^{h}\lambda^{g}\right)\cdot h, \eeq where \mat x^{h}_{g}\tam is defined
by the equation \beq{2} X(g)=\sum\limits_{h}x^{h}_{g}\cdot h\in\cA. \eeq Since
the sum in (\ref{2}) must be finite, this means that the matrix \mat
X=\|x^{h}_{g}\|_{g,h\in G}\tam must satisfy the following natural condition:
\begin{itemize}
\item[(\F1)] for every subscript $g\in G$, the set of superscripts $h\in G$ for
which $x^{h}_{g}$ is nonzero is finite.
\end{itemize}
In particular, it follows from condition (\F1) that, in the matrix
representation (\ref{1}), the outer sum is also finite.

Certainly, the converse also holds: if a matrix $X=\|x^{g}_{h}\|_{g,h\in G}$
satisfies condition (\F1), then it well defines a linear operator \mat
X\colon\cA\mapr{}\cA\tam by formula (\ref{1}). All this justifies the fact that
the operator \mat X\tam and its matrix \mat X=\|x^{h}_{g}\|_{g,h\in G}\tam are
denoted by the same symbol \mat X\tam.

Consider now the so-called derivation in the algebra $\cA$, i.e., an operator
$X$ satisfying the condition
\begin{itemize}
\item[(\F2)]$ X(u\cdot v)=X(u)\cdot v+u\cdot X(v), u,v\in \cA. $
\end{itemize}

The set of all derivations of the algebra $\cA$ is denoted by $\Der(\cA)$ and
forms a Lie algebra with respect to the commutator of operators.

A natural problem is to describe all derivations of the algebra $\cA$. To this
end, one should satisfy the conditions (\F1) and (\F2). The verification of
each of these conditions separately is a more or less simple task. The
simultaneous validity of these conditions is the content of the present paper.

There is a class of the so-called inner derivations, i.e., operators of the
form
$$
X=\ad(u), \quad X(v)=\ad(u)(v)=[u,v]=u\cdot v-v\cdot u, \quad  u,v\in \cA.
$$
All inner derivations satisfy both the conditions (\F1) and (\F2)
automatically. The set of these derivations is denoted by $\Int(\cA)$; it forms
a Lie subalgebra in the Lie algebra $\Der(\cA)$, \bdf
\Int(\cA)\subseteq\Der(\cA). \edf

\begin{proposition}
The Lie subalgebra \mat\Int(\cA)\subseteq\Der(\cA)\tam is an ideal.
\end{proposition}

Indeed, we are to verify the validity of the condition \bdf
[\Int(\cA),\Der(\cA)]\subset \Int(\cA). \edf If $\ad(u)\in\Int(\cA)$,
$X\in\Der(\cA)$, then the commutator $[\ad(u),X]$ is evaluated by the formula
 \bdf
\begin{array}{l}
[\ad(u),X](v)=\ad(u)(X(v))-X(\ad(u)(v))=[u,X(v)]-X([u,v])=
\\=[u,X(v)]-[X(u),v]-[u,X(v)]=-\ad(X(u))(v),
\end{array}
\edf which implies that $[\ad(u),X]\in\Int(\cA)$.

\section{Description of Derivations as Functions on the Groupoid $\cG$}

Denote by $\cG$ the groupoid associated with the adjoint action of the group
$G$ (or the corresponding action groupoid, see, e.g., Ershov(2012)
\cite{Ershov-2012}, p.~18, Example~j).

The groupoid $\cG$ consists of the objects $\Obj(\cG)=G$ and the morphisms \bdf
\Mor(a,b)=\{g\in G: ga=bg \hbox{ or } b=\Ad(g)(a)\}, \quad a,b\in \Obj(\cG).
\edf It is convenient to denote elements of the set of all morphisms
$\Mor(\cG)=\coprod\limits_{a,b\in \Obj(\cG)}\Mor(a,b)$ in the form of columns
\bdf \xi=\left(\frac{a\mapr{}b}{g}\right)\in \Mor(a,b), \quad
b=gag^{-1}=\Ad(g)(a). \edf The composition $*$ of two morphisms is defined by
the formula \bdf
\begin{array}{l}
\left(\frac{a\mapr{}c}{g_{2}g_{1}}\right)=
\left(\frac{b\mapr{}c}{g_{2}}\right)*
\left(\frac{a\mapr{}b}{g_{1}}\right),\\\\
b=\Ad(g_{1})(a), \\\\
c=\Ad(g_{2})(b)=\Ad(g_{2})(\Ad(g_{1})(a))=\Ad(g_{2}\Ad(g_{1})(a)),
\end{array}
\edf which corresponds to the diagram \bdf \xymatrix{
& \Ad(g_{1})(a)\ar@{=}[d]& \Ad(g_{2}g_{1})(a)\ar@{=}[d]\\
a\ar[r]^{g_{1}}\ar@/_15pt/[rr]_{g_{2}g_{1}}& b\ar[r]^{g_{2}}& c } \edf

There is another symbol for the morphism:
$$
\xi=\left(\xymatrix{a\ar[r]^{g}_{ga=bg}& b }\right)
$$
and for the composition of two morphisms: \bdf \xymatrix{
a\ar[rr]^{g_{1}}_{g_{1}a=bg_{1}}\ar@/^40pt/[rrrr]^{g_{2}g_{1}}_{g_{2}g_{1}a=cg_{2}g_{1}}&&
b\ar[rr]^{g_{2}}_{g_{2}b=cg_{2}}&& c } \edf

\subsubsection*{Operators as functions on the groupoid.}

A linear operator $X\colon\cA\mapr{}\cA$ is described by a matrix
$X=\|x^{h}_{g}\|_{g,h\in G}$ satisfying the following condition:
\begin{itemize}
\item[(\F1)] for every subscript $g\in G$, the set of all superscripts $h\in G$
for which $x^{h}_{g}$ is nonzero is finite.
\end{itemize}
The matrix $X=\|x^{h}_{g}\|_{g,h\in G}$ defines a function on the
groupoid~$\cG$,
$$
T^{X}\colon\Mor(\cG)\mapr{}R,
$$
associated with the operator~$X$; this function is defined by the following
formula: if \bdf \xi=\left(\frac{a\mapr{}b}{g}\right)\in\Mor(\cG), \edf then we
set \bdf T^{X}(\xi)=T^{X}\left(\frac{a\mapr{}b}{g}\right)=x^{ga=bg}_{g}. \edf
The condition (\F1) imposed on the coefficients of the matrix \mat X\tam can be
reformulated in terms of the function \mat T\tam:
\begin{itemize}
\item[(\T1)] for every element $g\in G$, the set of morphisms of the form \bdf
\xi=\left(\frac{a\mapr{}b}{g}\right) \edf for which \mat T^{X}(\xi)\neq 0\tam
is finite.
\end{itemize}
The set of all morphisms \mat\Mor(\cG)\tam can be represented in the form of a
disjoint union \bdf \Mor(\cG)=\coprod\limits_{g\in G}\cH_{g}, \edf where \bdf
\cH_{g}=\left\{\xi=\left(\frac{a\mapr{}b}{g}\right):a\in G,b=gag^{-1}\in
G\right\}. \edf

Then condition (\T1) imposed on the function \mat T\tam can be reformulated in
an equivalent way as follows.

\begin{proposition}
A function \bdf T^{X}\colon\Mor(\cG)\mapr{}\C \edf is defined by a linear
operator \bdf X\colon\cA\mapr{}\cA \edf if and only if, for every element \mat
g\in G,\tam the restriction \mat
{\left(T^{X}\right)_{|}}_{\cH_{g}}\colon\cH_{g}\mapr{}\C\tam is a finitely
supported function.
\end{proposition}

The functions \mat T\colon\Mor(\cG)\mapr{}\C\tam satisfying the finiteness
condition on every subset \mat\cH_{g},\tam \mat g\in G,\tam are called locally
finitely supported functions on the groupoid~\mat\cG\tam.

Consider two morphisms $\xi=\left(\frac{a\mapr{}b}{g}\right)$ and
$\eta=\left(\frac{b\mapr{}c}{g'}\right)$, which hence admit the composition
$$
\eta*\xi=\left(\frac{a\mapr{}c}{g'g}\right).
$$

\begin{theorem}
An operator $X\colon\cA\mapr{}\cA$ is a derivation if and only if the function
$T^{X}$ on the groupoid $\cG$ associated with the operator $X$ satisfies the
condition
\begin{itemize}
\item[\textup(\T2\textup)] $T^{X}(\eta*\xi)=T^{X}(\eta)+T^{X}(\xi)$
\end{itemize}
for every pair of morphisms $\xi$ and $\eta$ admitting the
composition~$\eta*\xi$.
\end{theorem}

\proof\ Let the matrix of the operator \mat X\tam be of the form \mat
X=\|x^{h}_{g}\|_{g,h\in G};\tam thus, the function \mat T^{X}\tam takes the
value \bdf T^{X}(\xi)=T^{X}\left(\frac{a\mapr{}b}{g}\right)=x^{ga=bg}_{g}. \edf
Let \mat \xi=\left(\frac{a\mapr{}b}{g_{1}}\right),\tam
\mat\eta=\left(\frac{b\mapr{}c}{g_{2}}\right),\tam \mat \eta*\xi=
\left(\frac{a\mapr{}c}{g_{2}g_{1}}\right).\tam Then \bdf
T^{X}(\eta*\xi)=x_{g_{2}g_{1}}^{g_{2}g_{1}a=cg_{2}g_{1}}=x_{g_{2}g_{1}}^{h},
\edf \bdf T^{X}(\xi)=x_{g_{1}}^{g_{1}a=bg_{1}}=x_{g_{1}}^{g_{2}^{-1}h}, \edf
\bdf T^{X}(\eta)=x_{g_{2}}^{g_{2}b=cg_{2}}=x_{g_{2}}^{hg_{1}^{-1}}. \edf On the
other hand, \bdf X(g_{2}g_{1})=X(g_{2})g_{1}+g_{2}X(g_{1}). \edf In other
words, \bdf
\begin{array}{l}
X(g_{2}g_{1})= \sum\limits_{h\in G} x^{h}_{g_{2}g_{1}}\cdot h=
\sum\limits_{h\in G}x^{h}_{g_{2}}\cdot h\cdot g_{1}+
g_{2}\cdot\sum\limits_{h\in G}x^{h}_{g_{1}}\cdot h\\= \sum\limits_{h\in
G}x^{hg_{1}^{-1}}_{g_{2}}\cdot h+ \sum\limits_{h\in
G}x^{g_{2}^{-1}h}_{g_{1}}\cdot h.
\end{array}
\edf Therefore, \bdf x^{h}_{g_{2}g_{1}}=
x^{hg_{1}^{-1}}_{g_{2}}+x^{g_{2}^{-1}h}_{g_{1}} \edf Thus, \bdf
T^{X}(\eta*\xi)=T^{X}(\eta)+T^{X}(\xi).\qed \edf

We refer to a function \mat T\colon\Mor(\cG)\mapr{}R\tam on the groupoid
\mat\cG\tam satisfying the additivity condition (T2) as a \textit{character};
denote the set of all characters on the groupoid \mat\cG\tam by
\mat\T(\cG).\tam Denote the space of all locally finitely supported characters
of the groupoid \mat\cG\tam by \mat\T_{f}(\cG)\subset\T(\cG).\tam

Thus, there is a mapping \bdf \Der(\cA)\mapr{T}\T_{f}(\cG), \edf which is
one-to-one.

\section{Inner derivations}
There are works (see, e.g., Losert(2008) \cite{Losert-2008}) related to the
so-called inner derivations of the group algebra. The commutator in the algebra
is a derivation, which is called an inner derivation.

A natural question arises: How are the inner derivations described in terms of
the matrix of the operator of the derivation?

The answer can be formulated as follows. Let \blf a\in G\elf, and let
\blf\ad(a)\elf be the commutator, \bdf \ad(a)(x)=[a,x],\quad x\in
C^{\infty}(G). \edf

This is an inner derivation. Denote by \mat \|A^{h}_{g}\|\tam the matrix of the
derivation \mat\ad(a).\tam Then \bdf \ad(a)(g)=\sum\limits_{h\in G}
A^{h}_{g}\cdot h \edf Since \mat\ad(a)(g)=ag-ga,\tam it follows that \bdf
A^{h}_{g}=\delta^{ag}_{h}-\delta^{ga}_{h}. \edf

The matrix of the operator \blf\ad(a)\elf defines the function \mat
T^{\ad(a)}\tam on the set of all morphisms \blf\Mor(\cG)\elf of the category
\blf\cG\elf. Let $\xi=\left(\frac{\alpha\mapr{}\beta}{g}\right)$ be a morphism
in the category \mat\cG\tam and let \mat g\alpha=\beta g(=h).\tam Then

\bdf T^{\ad(a)}(\xi)=T^{\ad(a)}\left(\frac{\alpha\mapr{}\beta}{g}\right)=
A^{g\alpha=\beta g}_{g}=\delta^{ag}_{g\alpha=\beta
g}-\delta^{ga}_{g\alpha=\beta g}. \edf

The first summand in the function \mat T^{\ad(a)}(\xi)\tam is equal to one if
and only if \blf \beta=a\elf, i.e., if and only if the morphism \blf\xi\in
\Mor(g^{-1}ag,a)\elf. Similarly, the second summand in the function \mat
T^{\ad(a)}(\xi)\tam is equal to minus one if and only if \blf \alpha=a\elf,
i.e., if and only if \blf\xi\in \Mor(a,gag^{-1})\elf.

In other words, the matrix \blf\ad(a)\elf is equal to one on the
morphisms\break \blf \Mor(g^{-1}ag,a),\elf to minus one on the morphisms \elf
\Mor(a,gag^{-1}),\elf and is equal to zero on the morphisms \elf \Mor(u,u)\elf
and on the morphisms \elf \Mor(a,a)\elf and \elf \Mor(v,v)\elf in the following
diagram:

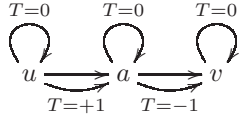
\bdf \xymatrix{ u\ar@(lu,ru)[]^{T=0}\ar[r]\ar@/_/[r]_{T=+1}
&a\ar@/_/[r]_{T=-1}\ar[r] \ar@(lu,ru)[]^{T=0}&v\ar@(lu,ru)[]^{T=0} } \edf

This proves the following theorem.

\begin{theorem}[on the inner derivations]
\end{theorem}
The characters of the inner derivations are trivial on \mat\Mor(a,a)\tam: \bdf
\xymatrix{
\Int(\cA)\ar[r]\ar[d]_{\cap}&\Der(\cA)\ar[d]^{T}_{\approx}\\
\ker p_{a}\ar[r]& \T_{f}(\cG)\ar[r]^{p_{a}}& \T_{f}(\Mor(a,a)) } \edf The set
\mat\T_{f}(\Mor(a,a))\tam coincides with the group of all characters \bdf
\T_{f}(\Mor(a,a))=\T(\Mor(a,a)). \edf

Note that, if a character \mat T\in \T(\cG)\tam vanishes on \mat\Mor(a,a),\tam
then it vanishes on \mat\Mor(u,u)\tam for every conjugate element \mat u\in[a],
\quad u=gag^{-1}.\mat

Thus, the diagram has the following form: \bdf \xymatrix{
\Int(\cA)\ar[r]\ar[d]_{\cap}^{T}&\Der(\cA)\ar[d]^{T}_{\approx}\\
\ker p_{a}\ar[r]& \T_{f}(\cG)\ar[r]^{p_{a}}& \T(\Mor(a,a)) } \edf

From the viewpoint of the Johnson derivation problem
(Johnson(2001)\cite{Johnson-2001}), it is natural to denote by
\mat\Out(\cA)\tam the quotient group \mat\Out(\cA)=\Der(\cA)/\Int(\cA)\tam and
refer to it as the algebra of outer derivations of the algebra \mat\cA\tam.
Thus, the previous diagram is completed to the diagram \bdf \xymatrix{
0\ar[r]&\Int(\cA)\ar[r]\ar[d]_{\cap}^{T}&
\Der(\cA)\ar[d]^{T}_{\approx}\ar[r]&\Out(\cA)\ar[d]\ar[r]&0\\
0\ar[r]&\ker p_{a}\ar[r]& \T_{f}(\cG)\ar[r]^{p_{a}}& \T(\Mor(a,a)) } \edf

\subsection*{Description of inner derivations.}

Let us note first that the set of morphisms \mat\Mor(\cG)\tam of the groupoid
\mat\cG\tam is decomposed into a disjoint union of morphisms over the conjugacy
classes of \mat G\tam which are the objects of the groupoid \mat\cG\tam by
definition. The group \mat G\tam is decomposed into the disjoint union of the
conjugacy classes \bdf G=\coprod\limits_{g\in G}[g], \quad [g]=\{h: \exists
a\in G, h=aga^{-1}\}. \edf Correspondingly, the set of morphisms is also
represented as the disjoint union \bdf
\Mor(\cG)=\coprod\limits_{[g]}\Mor(\cG_{[g]}). \edf This means that the
construction of every derivation can be carried out in the form of derivations
\mat\Der_{[g]}(\A)\tam independently in every subcategory \mat\cG_{[g]}\mat \
as locally finitely supported characters on each of these subcategories.

A natural problem is to find out whether or not the set of all derivations
trivial on all \blf\Mor(u,u)\elf coincides with the set of inner derivations.
In other words, whether or not the embedding \mat\xymatrix{ \Int_{\hskip
-4pt[u]}(\cA)\ar[r]^{\subset}& \ker p_{u}}\tam is an isomorphism:
\bdf\xymatrix{ \Int_{\hskip -4pt [u]}(\cA)\ar[r]^{\subset}& \ker
p_{u}\ar[r]^{\subset} &\T_{f}(\cG)\ar[r]^{p_{u}}& \T(\Mor(u,u)) }. \edf

The investigation of this problem enables us to formulate specific conditions
on a locally finitely supported character \mat T\colon\Mor(\cG)\mapr{}R\tam
that realizes a given inner derivation \mat X\in\Int(\cG)\tam, \bdf T=T^{X}.
\edf

\subsubsection*{Case of the identity element \mat[e].\tam}

In particular, one of the subcategories corresponds to the identity element
\mat e\in G\tam for which \mat[e]=\{e\}\tam. In this special case, the
subcategory \mat\cG_{[e]}\tam consists of a single object \mat e\in G\tam and
the set of morphisms is isomorphic to the group \mat G,\tam
\mat\Mor(e,e)\approx G.\tam In particular, the set of locally finitely
supported characters \mat T_{f}(\cG_{[e]})\approx T(\cG_{[e]})\approx T(G)\tam
is isomorphic to the group of all characters on the group \mat G.\tam Every
character on the group, \mat T\in T(G)\tam, is realized as a derivation \mat
X\in\Der(G),\tam \mat T^{X}=T\tam. Indeed, the character \mat T\in T(G)\tam is
a character on the category \mat\cG\tam which is equal to \mat T\tam on
\mat\Mor(e,e)\tam and to zero on the other summands \mat\Mor(\cG_{[g]}), \quad
g\neq e\tam. Therefore, the corresponding matrix \mat \| X^{h}_{g}\|\tam of the
operator \mat X\tam is given by the formula \bdf X^{h}_{g}=T(g)\delta^{h}_{g}.
\edf All derivations corresponding to the characters on the subcategory
\mat\cG_{[e]}\tam are not inner derivations.

Similar considerations fit for the other conjugacy classes that consist of
finitely many elements, i.e., when \mat\#[g]<+\infty.\mat \ In particular, this
holds for the elements in the center \mat g\in Z(G).\tam

\subsubsection*{Exact sequence.}

If we get rid of the condition that the characters are locally finitely
supported, then one can establish that some sequence is exact, as is formulated
in the following theorem.

\begin{theorem}\label{t3}
The following sequence is exact\/\textup: \bdf 0\mapr{}\ker
p_{a}\mapr{}\T(\cG_{[a]})\mapr{p_{a}}\T(\Mor(a,a))\mapr{}0, \quad a\in U. \edf
\end{theorem}

\proof\ One should prove only that the mapping \blf p(a)\elf is epimorphic.

Let \blf\chi\in\X(G^{a})\elf.

Choose an element \blf a\in U\elf and arbitrary elements \blf
g_{a,b}\in\Mor(a,b),\forall b\in U \elf such that \blf g_{a,a}=e\elf. Elements
of this kind exist indeed, because \blf\Mor(a,b)\elf is nonempty. Write further
\bdf g_{b,c}=g_{a,c}g_{a,b}^{-1}\in\Mor(b,c). \edf Note that the condition \blf
g_{b,c}\in\Mor(b,c)\elf means that \bdf g_{b,c}b=c g_{b,c}. \edf Let us prove
this condition. First, \bdf g_{a,b}a=b g_{a,b}\edf by construction.  Further,
\blf g_{b,a}=g_{a,a}g_{a,b}^{-1}=g_{a,b}^{-1}\in\Mor(b,a),\elf i.e., \bdf
g_{b,c}=g_{a,c}g_{b,a}. \edf Thus, \bdf
g_{b,c}b=g_{a,c}g_{b,a}b=g_{a,c}ag_{b,a}= c g_{a,c}g_{b,a}=c g_{b,c}, \edf
i.e., \bdf g_{b,c}\in\Mor(b,c). \edf

Hence, the following relation holds: \bdf g_{c,d}g_{b,c}=g_{b,d}\in\Mor(b,d).
\edf

Indeed, \bdf g_{c,d}g_{b,c}=g_{a,d}g_{c,a}g_{a,c}g_{b,a}=
g_{a,d}g_{b,a}=g_{b,d}\in\Mor(b,d). \edf

Let us now construct a character \blf X\elf on the groupoid \blf\cG_{U},\elf
\blf X\colon\Mor(U)\mapr{}\R\elf. Let \blf x_{b,c}\in\Mor(b,c)\elf be an
arbitrary morphism. Then \bdf g_{c,a}x_{b,c}g_{a,b}\in\Mor(a,a)=G^{a}. \edf
Write \bdf X(x_{b,c})\bydef\chi(g_{c,a}x_{b,c}g_{a,b}). \edf It can readily be
seen that the mapping \blf X\colon\Mor(U)\mapr{}\R\elf is additive, \bdf
\begin{array}{l}
X(x_{c,d}x_{b,c})=\chi(g_{d,a}(x_{c,d}x_{b,c})g_{a,b})=
\chi(g_{d,a}x_{c,d}g_{a,c}g_{c,a}x_{b,c}g_{a,b})\\=
\chi(g_{d,a}x_{c,d}g_{a,c})\chi(g_{c,a}x_{b,c}g_{a,b})= X(x_{c,d})X(x_{b,c}).
\end{array}
\edf

The restriction of the character \blf X\elf to \blf\Mor(a,a)=G^{a}\elf
coincides with \blf\chi\elf, \bdf
X(x_{a,a})=\chi(g_{a,a}x_{a,a}g_{a,a})=\chi(x_{a,a}).\qed \edf

\subsubsection*{Reduction to Groups of Cochains}
Let us return to the study of characters \mat T\tam that are trivial on the
subspace \mat\Mor(a,a)\tam, \bdf\xymatrix{ \Int_{a}\cA\ar[r]^{\subset}& \ker
p^{f}_{a}\ar[r]^{\subset} &\T_{f}(\cG_{[a]})\ar[r]^{p^{f}_{a}}& \T(\Mor(a,a)).
} \edf Denote by \mat\Delta(\cG_{[a]})\tam the simplex whose vertices are the
elements of the conjugacy class \mat[a]\subset G=\Obj(\cG_{[a]}).\tam

Since every character \mat T\in\ker p^{f}_{a}\tam takes equal values on the set
of all morphisms \mat\Mor(b,c)\tam, \mat b,c\in[a]\subset G\tam, it follows
that there is a natural embedding \bdf \varphi\colon\ker
p^{f}_{a}\hookrightarrow \C^{1}(\Delta(\cG_{[a]}))\edf in the group of cochains
of the simplex \mat\Delta(\cG_{[a]})\tam, and every character \mat T\in\ker
p^{f}_{a}\tam is a cocycle: \bdf \xymatrix{
\C^{0}_{f}(\Delta(\cG_{[a]}))\ar@{^{(}->}[r]&\C^{0}(\Delta(\cG_{[a]}))\ar[r]^{\delta}&\C^{1}(\Delta(\cG_{[a]}))\ar[r]^{\delta}&C^{2}(\Delta(\cG_{[a]}))\\
\Int_{a}(\cA)\ar@{^{(}->}[r]^{T}\ar@{^{(}->}[ur]\ar@{^{(}->}[u]^{\varphi_{i0}}&
\ker
p^{f}_{a}\ar@{^{(}->}[u]^{\varphi_{0}}\ar@{^{(}->}[ur]^{\varphi}\ar[r]&0\ar[ur]
} \edf

The embedding \mat \varphi_{i0}\tam is an isomorphism. The image \bdf
\delta(\varphi_{i0}(\Int_{a}(\cA)))\subset
\delta(C^{0}_{f}(\Delta(\cG_{[a]})))\subset \ker p^{f}_{a}\subset
\C^{1}(\Delta(\cG_{[a]})) \edf can be described as some set of cocycles
satisfying certain conditions. Consider the conjugacy class \mat[a]\tam on
which the group \mat G\tam acts by the adjoint action \bdf
G\times[a]\mapr{\ad}[a], \quad \ad_{g}(b)=gbg^{-1}, \, b\in[a]. \edf For an
arbitrary element \mat g\in G\tam, consider the graph \mat\Gamma_{g,a}\subset
\Delta(\cG_{[a]}),\tam formed by directed edges whose beginnings are the
elements \mat b\in[a]\tam and the ends are \mat gbg^{-1}\in[a]\tam. Thus, every
edge is of the form of a directed segment \mat[b, gbg^{-1}].\tam The graph
\mat\Gamma_{g,a}\tam is decomposed into a disjoint union of directed paths formed
by the directed edges \bdf
\Gamma_{g,a}=\coprod\limits_{\alpha}\Gamma_{g,a}^{\alpha}. \edf Every directed
path \mat\Gamma_{g,a}^{\alpha}\tam can be infinite in both the directions or
finite, in which case this path is cyclic.

\begin{theorem}
Let \mat X\in\Int_{[a]}(\cA).\tam Then the cochain \mat
\varphi(T^{X})\in\C^{1}(\Delta(\cG_{[a]}))\tam satisfies the condition
\begin{itemize}
\item[\F\F:] the cochain \mat \varphi(T^{X})\tam is finitely supported on the
graph \mat\Gamma_{g,a},\tam and the sum of values of the cochain \mat
\varphi(T^{X})\tam on every directed path \mat\Gamma_{g,a}^{\alpha}\tam
vanishes.
\end{itemize}\end{theorem}

It is not clear whether the converse is true that is
if an operator \mat X\in\Der(\cA),\tam \mat T^{X}\in\ker
p^{f}_{a}\subset T_{f}(\cG_{[a]})\tam, satisfies condition (\F\F), whether \mat
X\in\Int_{[a]}(\cA)? \tam

\subsubsection*{The kernel \mat\ker p^{f}_{a}\tam differs from the set of inner derivations}

There are examples of groups for which there are locally finitely supported
characters \mat T\in \ker p^{f}_{a}\tam that do not satisfy condition (\F\F).
For the simplest example, consider the free group with two generators \mat
G=F<x_{1},x_{2}>.\tam For the conjugacy class we take the class
\mat[x_{1}]\subset G.\tam Consider the character \mat
T:\Mor(\cG_{[x_{1}]})\mapr{} R\tam, defining the values of the character on the
generators of the groupoid \mat\Mor(\cG_{[x_{1}]})\tam independently of one
another. The set of generators of the groupoid \mat\Mor(\cG_{[x_{1}]})\tam
consists of the morphisms of the form  \bdf
\xi=\left(\frac{\alpha\mapr{}\beta}{g}\right), \quad g=x_{1},x_{2}, \quad
\beta=g\alpha g^{-1}, \quad \alpha\in[x_{1}]. \edf Set \bdf
\begin{array}{l}
T\left(\frac{x_{2}x_{1}x_{2}^{-1}\mapr{}x_{1}x_{2}x_{1}x_{2}^{-1}x_{1}^{1}}{x_{1}}\right)=1;\\
\hbox{and let the character }T\hbox{ vanish on the other generators}.
\end{array}
\edf Since in the free group there are no relations except for natural
reductions in words, it follows that the function \mat T\tam can be extended by
additivity to some character on the groupoid \mat\cG_{[x_{1}]}.\tam

Hence, the character \mat T\tam takes only one value equal to 1 on one of the
directed paths of the form \mat\Gamma^{\alpha}_{x_{1},x_{1}},\tam and this
character is identically equal to zero on the other directed paths. This means
that condition (\T1) holds, and condition (\F\F) fails to hold.

\section{Appendix: Groups with Finitely Many Generators}

Let a group $G$ be finitely presentable, i.e., let it have finitely many
generators $\{g_{1}, g_{2}, \dots, g_{k}\}$ and finitely many defining
relations $\{S_{1}, S_{2}, \dots, S_{l}\}$:
$$
G=F<g_{1}, g_{2}, \dots, g_{k}>/\{S_{1}, S_{2}, \dots, S_{l}\}.
$$

Every relation \mat S_{i}\tam is a word of length \mat s_{i}\tam formed of
generating elements \mat g_{1}, g_{2}, \dots, g_{k}\tam or inverses of the
generating elements \mat g_{1}^{-1}, g_{2}^{-1}, \dots, g_{k}^{-1}\tam. Thus,
the relation \mat S_{i}\tam can be represented in the following form: \bdf
S_{i}=g_{j_{(i,1)}}^{\varepsilon_{(i,1)}} g_{j_{(i,2)}}^{\varepsilon_{(i,2)}}
g_{j_{(i,3)}}^{\varepsilon_{(i,3)}} \dots
g_{j_{(i,s_{i})}}^{\varepsilon_{(i,s_{i})}}\equiv 1, \quad
\varepsilon_{(i,j)}=\pm 1. \edf

Every relation \mat S_{i}\tam induces a series of relations on the groupoid
\mat\cG\tam: \bdf \xi_{(i,1)}\xi_{(i,1)}\xi_{(i,3)}\dots\xi_{(i,s_{i})}\equiv
\left(\frac{\alpha_{(i,1)}\mapr{}\alpha_{(i,1)}} {1}\right), \edf where the
morphisms \mat \xi_{(i,j)}\tam are defined by the rule \bdf
\begin{array}{l}
\xi_{(i,j)}=\left(\frac{\alpha_{(i,j)}\mapr{}\beta_{(i,j)}}
{g_{j_{(i,j)}}^{\varepsilon_{(i,j)}}}\right),\\\\
\beta_{(i,j)}=g_{j_{(i,j)}}^{\varepsilon_{(i,j)}}\alpha_{(i,j)}
g_{j_{(i,j)}}^{-\varepsilon_{(i,j)}},\\\\
\beta_{(i,j)}=\alpha_{(i,j+1)},\quad \beta_{(i,s_{i})}=\alpha_{(i,1)}.
\end{array}
\edf

Thus, to define a locally finitely supported character \mat T\tam on the
groupoid \mat\cG,\tam it suffices to define the values of the locally finitely
supported character \mat T\tam on the set of generators \bdf
\coprod\limits_{i=1}^{k}\cH_{g_{i}} \edf in such a way that the additivity
condition  \bdf
T(\xi_{(i,1)})+T(\xi_{(i,1)})+T(\xi_{(i,3)}+\cdots+T(\xi_{(i,s_{i})})=0 \edf
holds on every relation of the form \bdf
\xi_{(i,1)}\xi_{(i,1)}\xi_{(i,3)}\dots\xi_{(i,s_{i})}\equiv
\left(\frac{\alpha_{(i,1)}\mapr{}\alpha_{(i,1)}} {1}\right). \edf

\begin{otherlanguage}{russian}

\newtheorem{theoremru}{Теорема}
\newtheorem{propositionru}{Предложение}
\newtheorem{corru}{Следствие}
\def\proofru{{\bf Докзательство.}}

\title{
Деривации групповых алгебр}
\author{Арутюнов, А.А.\thanks{Работа выполнена при финансовой поддержке МОН РФ (Соглашение № 02.a03.21.0008  от 24.06.2016).}, Мищенко, А.С.\thanks{частично поддержан грантом РФФИ № 14-01-00007}, Штерн, А.И.\thanks{частично поддержан грантом РФФИ № 14-01-00007}}

\begin{dedication}{Памяти Ю.П.Соловьева посвящается}\end{dedication}
\maketitle

\begin{abstract}
В работе дается метод описания внешних дериваций групповой алгебры конечно представимой группы. Описание дериваций дается в терминах характеров группоида присоединенного действия группы.
\end{abstract}

\section{Введение}

Рассмотрим алгебру \mat\cA\tam  и некоторый бимодуль \mat E\tam над алгеброй \mat\cA\tam. Обозначим через \mat\Der(\cA,E)\tam пространство всех дериваций из алгебры \mat\cA\tam в бимодуль \mat E\tam, т.е. множество отображений
\bdf
D:\cA\mapr{}E,
\edf
которые удовлетворяют условию:
\bdf
D(ab)=D(a)b+aD(b), \quad a,b\in\cA,
\edf
(см. Losert(2008) \cite{Losert-2008}, Ghahramani(2000) \cite{Ghahramani-2000}). Среди дериваций \mat\Der(\cA,E)\tam выделяются так называемые внутренние деривации
\mat\Int(\cA,E)\subset\Der(\cA,E) \tam, которые задаются присоединенными представлениями
\bdf
\ad_{x}(a)\bydef xa-ax, \quad x\in E, a\in\cA.
\edf

Проблема дериваций формулируется следующим образом: все ли деривации являются внутренними? Эта задача рассматривалась не для всяких алгебр, а для групповых алгебр \mat\cA=C[G]\tam некоторой группы \mat G\tam. Более точно, рассматривается групповая алгебра \mat \overline\cA=L^{1}(G)\tam и бимодуль \mat E=M(G),\tam где \mat M(G)\tam есть алгебра всех ограниченных мер на группе \mat G\tam с операцией умножения, задаваемой сверткой мер.

Вопрос из работы Dales(2000) \cite{Dales-2000}, (Question  5.6.B, стр.746) формулируется следующим образом: Пусть \mat G\tam есть локально компактная группа. Всякая ли деривация из алгебры \mat \cA=L^{1}(G)\tam в бимодуль \mat E=M(G)\tam является внутренней деривацией? Утвердительный ответ оправдывается следующим соображением.

В случае, когда группа \mat G\tam является дискретной свободной абелевой группой с конечным числом образующих, т.е. \mat G\approx\mathbb{Z}^{n}\mat, то алгебру
\mat \overline\cA=L^{1}(G)\tam можно можно отождествить с алгеброй Фурье \mat A(\mathbb{T}^{n})\tam непрерывных функций на \mat n\tam--мерном торе \mat\mathbb{T}^{n},\tam коэффициенты Фурье которых образуют абсолютно сходящийся кратный ряд, \mat \cA =A(\mathbb{T}^{n})\subset C(\mathbb{T}^{n}),\tam  (эта алгебра Фурье  меньше алгебры непрерывных функций). Дериваций на алгебре \mat A(\mathbb{T}^{n})\tam нет, поскольку в ней достаточно много негладких функций, впрочем и внутренних дериваций тоже нет, поскольку алгебра \mat \overline\cA=L^{1}(G)\tam коммутативна.

Нас же интересует не вся банахова  алгебра \mat \overline\cA=L^{1}(G),\tam а только ее плотная подалгебра \mat \cA=C[G]\subset\overline\cA\mat, состоящая, так сказать, из гладких элементов в алгебре \mat \overline\cA=L^{1}(G),\tam следуя терминологии А.Кона (\cite{Connes-1994}, стр. 247). Для групповой алгебры \mat \cA=C[G]\;\mat тоже можно сформулировать аналогичную задачу: описать алгебру всех внешних дериваций групповой алгебры \mat \cA=C[G].\mat

\section{Групповая алгебра $C[G]$.}

Рассмотрим групповую алгебру
$\cA=C[G].$ Мы предполагаем, что группа \mat G\tam является конечно представимой дискретной группой.

Произвольный элемент $u\in\cA$ --- это конечная линейная комбинация
$
u=\sum\limits_{g\in G}\lambda^{g}\cdot g.
$
Рассмотрим произвольный линейный оператор
$$
X:\cA\mapr{}\cA
$$
Линейный оператор $X$ имеет следующий матричный вид
\beq{1}
X(u)=\sum\limits_{h\in G}\left(\sum\limits_{g\in G}x_{g}^{h}\lambda^{g}\right)\cdot h,
\eeq
где \mat x^{h}_{g}\tam определяется равенством
\beq{2}
X(g)=\sum\limits_{h}x^{h}_{g}\cdot h\in\cA.
\eeq
Поскольку сумма в равенстве (\ref{2}) должна быть конечной, то это значит, что
матрица \mat X=\|x^{h}_{g}\|_{g,h\in G}\tam должна удовлетворять естественному условию:
\begin{itemize}
\item[(\F1)] Для любого индекса $g\in G$ множество тех индексов $h\in G$, для которых $x^{h}_{g}$ отлично от нуля, конечно.
\end{itemize}
В частности, из условия (\F1) следует, что в матричном представлении (\ref{1}) внешняя сумма тоже конечна.

Разумеется, верно и обратное утверждение: если матрица $X=\|x^{g}_{h}\|_{g,h\in G}$ удовлетворяет условию (\F1), то она корректно задает линейный оператор \mat X:\cA\mapr{}\cA\tam по формуле (\ref{1}). Все это оправдывает, что и оператор \mat X,\tam и его матрица \mat X=\|x^{h}_{g}\|_{g,h\in G}\tam обозначаются тем
же самым символом \mat X\tam.

Рассмотрим теперь так называемое дифференцирование (деривацию) в алгебре $\cA$, т.е. такой оператор
$X$, для которого выполнено условие
\begin{itemize}
\item[(\F2)]$
X(u\cdot v)=X(u)\cdot v+u\cdot X(v), u,v\in \cA.
$
\end{itemize}

Множество всех дериваций алгебры $\cA$ обозначается через $\Der(\cA)$ и образует
алгебру Ли по отношению к коммутатору операторов.

Естественная задача заключается в том, чтобы описать все дифференцирования алгебры $\cA$. Для этого нужно соблюсти два условия (\F1) и (\F2). Каждое условие в отдельности проверяется более или менее просто. Одновременное выполнение этих условий составляет содержание настоящей работы.

Имеется класс так называемых внутренних дифференцирований, т.е. операторов вида
$$
X=\ad(u), \quad X(v)=\ad(u)(v)=[u,v]=u\cdot v-v\cdot u, \quad  u,v\in \cA
$$
Все внутренние деривации автоматически удовлетворяют обоим условиям (\F1) и (\F2). Они обозначаются через $\Int(\cA)$ и образуют подалгебру Ли в алгебре Ли $\Der(\cA)$,
\bdf
\Int(\cA)\subseteq\Der(\cA).
\edf

\begin{propositionru}
Подалгебра \mat\Int(\cA)\subseteq\Der(\cA)\tam идеалом.
\end{propositionru}

Действительно, требуется проверить выполнение условия
\bdf
[\Int(\cA),\Der(\cA)]\subset \Int(\cA).
\edf
Если $\ad(u)\in\Int(\cA)$, $X\in\Der(\cA)$, то коммутатор $[\ad(u),X]$ вычисляется
по формуле:
\bdf
\begin{array}{l}
[\ad(u),X](v)=\ad(u)(X(v))-X(\ad(u)(v))=[u,X(v)]-X([u,v])=
\\=[u,X(v)]-[X(u),v]-[u,X(v)]=-\ad(X(u))(v),
\end{array}
\edf
т.е. $[\ad(u),X]\in\Int(\cA)$.

\section{Описание дериваций как функций на группоиде $\cG$}

Обозначим через $\cG$ группоид, ассоциированный с присоединенным действием группы $G$ (или группоид действия, см. например, Ершов(2012)
\cite{Ershov-2012}, стр. 18, пример j).

Группоид $\cG$ состоит из объектов $\Obj(\cG)=G$ и морфизмов
\bdf
\Mor(a,b)=\{g\in G: ga=bg \hbox{ или } b=\Ad(g)(a)\}, \quad a,b\in \Obj(\cG).
\edf
Элементы множества всех морфизмов
$\Mor(\cG)=\coprod\limits_{a,b\in \Obj(\cG)}\Mor(a,b)$ удобно обозначать в виде столбца
\bdf
\xi=\left(\frac{a\mapr{}b}{g}\right)\in \Mor(a,b), \quad b=gag^{-1}=\Ad(g)(a).
\edf
Композиция $*$ двух морфизмов задается формулой
\bdf
\begin{array}{l}
\left(\frac{a\mapr{}c}{g_{2}g_{1}}\right)=
\left(\frac{b\mapr{}c}{g_{2}}\right)*
\left(\frac{a\mapr{}b}{g_{1}}\right),\\\\
b=\Ad(g_{1})(a), \\\\
c=\Ad(g_{2})(b)=\Ad(g_{2})(\Ad(g_{1})(a))=\Ad(g_{2}\Ad(g_{1})(a))
\end{array}
\edf
которая соответствует диаграмме
\bdf
\xymatrix{
& \Ad(g_{1})(a)\ar@{=}[d]& \Ad(g_{2}g_{1})(a)\ar@{=}[d]\\
a\ar[r]^{g_{1}}\ar@/_15pt/[rr]_{g_{2}g_{1}}&
b\ar[r]^{g_{2}}&
c
}
\edf

Другое изображение морфизма:
$$
\xi=\left(\xymatrix{a\ar[r]^{g}_{ga=bg}& b }\right)
$$
и композиции двух морфизмов
\bdf
\xymatrix{
a\ar[rr]^{g_{1}}_{g_{1}a=bg_{1}}\ar@/^40pt/[rrrr]^{g_{2}g_{1}}_{g_{2}g_{1}a=cg_{2}g_{1}}&&
b\ar[rr]^{g_{2}}_{g_{2}b=cg_{2}}&&
c
}
\edf

\subsubsection*{Операторы как функции на группоиде.}

Линейный оператор $X:\cA\mapr{}\cA$
описывается матрицей
$X=\|x^{h}_{g}\|_{g,h\in G}$, которая  удовлетворяет условию:
\begin{itemize}
\item[(\F1)] Для любого индекса $g\in G$ множество  тех индексов $h\in G$, для которых $x^{h}_{g}$ отлично от нуля, конечно.
\end{itemize}
Матрица $X=\|x^{h}_{g}\|_{g,h\in G}$ задает функцию на группоиде $\cG$
$$
T^{X}:\Mor(\cG)\mapr{}R,
$$
ассоциированную с оператором $X$,
которая определяется формулой: если
\bdf
\xi=\left(\frac{a\mapr{}b}{g}\right)\in\Mor(\cG),
\edf
то полагаем
\bdf
T^{X}(\xi)=T^{X}\left(\frac{a\mapr{}b}{g}\right)=x^{ga=bg}_{g}.
\edf
Условие (\F1), налагаемое на коэффициенты матрицы \mat X\tam, можно переформулировать в терминах функции \mat T\tam:
\begin{itemize}
\item[(\T1)] Для любого индекса $g\in G$ множество морфизмов вида
\bdf
\xi=\left(\frac{a\mapr{}b}{g}\right),
\edf
для которых \mat T^{X}(\xi)\neq 0,\tam конечно.
\end{itemize}
Множество всех морфизмов \mat\Mor(\cG)\tam представляется в виде несвязного объединения
\bdf
\Mor(\cG)=\coprod\limits_{g\in G}\cH_{g},
\edf
где
\bdf
\cH_{g}=\left\{\xi=\left(\frac{a\mapr{}b}{g}\right):a\in G,b=gag^{-1}\in G\right\}.
\edf

Тогда условие (\T1), налагаемое на функцию \mat T\tam можно эквивалентным образом переформулировать следующим образом
\begin{propositionru}
Функция
\bdf
T^{X}:\Mor(\cG)\mapr{}\C,
\edf
задается линейным оператором
\bdf
X:\cA\mapr{}\cA
\edf
тогда и только тогда, когда для любого элемента \mat g\in G\tam
ограничение \mat {\left(T^{X}\right)_{|}}_{\cH_{g}}:\cH_{g}\mapr{}\C\tam является финитной функцией.
\end{propositionru}

Такие функции \mat T:\Mor(\cG)\mapr{}\C\tam, которые удовлетворяют условию финитности на каждом подмножестве \mat\cH_{g},\tam \mat g\in G\tam, будем называть локально финитными функциями на группоиде \mat\cG\tam.

Рассмотрим два морфизма $\xi=\left(\frac{a\mapr{}b}{g}\right)$ и
$\eta=\left(\frac{b\mapr{}c}{g'}\right)$ которые допускают, следовательно,
композицию
$$
\eta*\xi=\left(\frac{a\mapr{}c}{g'g}\right).
$$
\begin{theoremru}
Оператор $X:\cA\mapr{}\cA$ является дифференцированием (т.е. деривацией) тогда и только тогда, когда для ассоциированной с оператором $X$ функции $T^{X}$ на группоиде $\cG$ выполняется условие
\begin{itemize}
\item[(\T2)] $T^{X}(\eta*\xi)=T^{X}(\eta)+T^{X}(\xi)$
\end{itemize}
для любой пары морфизмов $\xi$ и $\eta$, допускающих композицию $\eta*\xi$.
\end{theoremru}
\proofru Пусть матрица оператора \mat X\tam имеет вид
\mat X=\|x^{h}_{g}\|_{g,h\in G},\tam значит, функция \mat T^{X}\tam принимает значение
\bdf
T^{X}(\xi)=T^{X}\left(\frac{a\mapr{}b}{g}\right)=x^{ga=bg}_{g}.
\edf
Пусть \mat \xi=\left(\frac{a\mapr{}b}{g_{1}}\right),\tam \mat\eta=\left(\frac{b\mapr{}c}{g_{2}}\right),\tam
\mat \eta*\xi= \left(\frac{a\mapr{}c}{g_{2}g_{1}}\right).\tam Тогда
\bdf
T^{X}(\eta*\xi)=x_{g_{2}g_{1}}^{g_{2}g_{1}a=cg_{2}g_{1}}=x_{g_{2}g_{1}}^{h},
\edf
\bdf
T^{X}(\xi)=x_{g_{1}}^{g_{1}a=bg_{1}}=x_{g_{1}}^{g_{2}^{-1}h},
\edf
\bdf
T^{X}(\eta)=x_{g_{2}}^{g_{2}b=cg_{2}}=x_{g_{2}}^{hg_{1}^{-1}}.
\edf
С другой стороны,
\bdf
X(g_{2}g_{1})=X(g_{2})g_{1}+g_{2}X(g_{1}).
\edf
Другими словами,
\bdf
\begin{array}{l}
X(g_{2}g_{1})=
\sum\limits_{h\in G} x^{h}_{g_{2}g_{1}}\cdot h=
\sum\limits_{h\in G}x^{h}_{g_{2}}\cdot h\cdot g_{1}+
g_{2}\cdot\sum\limits_{h\in G}x^{h}_{g_{1}}\cdot h=\\\\=
\sum\limits_{h\in G}x^{hg_{1}^{-1}}_{g_{2}}\cdot h+
\sum\limits_{h\in G}x^{g_{2}^{-1}h}_{g_{1}}\cdot h.
\end{array}
\edf
Значит,
\bdf
x^{h}_{g_{2}g_{1}}=
x^{hg_{1}^{-1}}_{g_{2}}+x^{g_{2}^{-1}h}_{g_{1}}
\edf
Таким образом,
\bdf
T^{X}(\eta*\xi)=T^{X}(\eta)+T^{X}(\xi).
\edf
\qed

Функцию \mat T:\Mor(\cG)\mapr{}R\tam на группоиде \mat\cG,\tam удовлетворяющую условию аддитивности (T2), будем называть характером,
а множество всех характеров на группоиде \mat\cG\tam обозначать через
\mat\T(\cG).\tam Пространство всех  локально финитных характеров группоида \mat\cG\tam будем обозначать через \mat\T_{f}(\cG)\subset\T(\cG).\tam

Таким образом существует отображение
\bdf
\Der(\cA)\mapr{T}\T_{f}(\cG),
\edf
которое является взаимно однозначным отображением.

\section{О внутренних деривациях}
Имеются некоторые работы (см., например, Losert(2008) \cite{Losert-2008}), связанные с так называемым внутренними дифференцированиями или внутренними деривациями групповой алгебры. Коммутатор в алгебре является деривацией, которая называется внутренней деривацией.

Возникает естественный вопрос: как внутренние деривации описываются в терминах матрицы оператора деривации?

Ответ можно сформулировать следующим образом. Пусть \blf a\in G\elf, \blf\ad(a)\elf -- коммутатор:
\bdf
\ad(a)(x)=[a,x],\quad x\in C^{\infty}(G).
\edf

Это внутренняя деривация.
Обозначим через \mat \|A^{h}_{g}\|\tam матрицу деривации
\mat\ad(a).\tam Тогда
\bdf
\ad(a)(g)=\sum\limits_{h\in G} A^{h}_{g}\cdot h
\edf
Поскольку \mat\ad(a)(g)=ag-ga,\tam то
\bdf
A^{h}_{g}=\delta^{ag}_{h}-\delta^{ga}_{h}.
\edf

Матрица оператора \blf\ad(a)\elf задает функцию \mat T^{\ad(a)}\tam  на множестве всех морфизмов \blf\Mor(\cG)\elf категории \blf\cG\elf. Пусть
$\xi=\left(\frac{\alpha\mapr{}\beta}{g}\right)$ -- морфизм в категории \mat\cG,\tam
\mat g\alpha=\beta g(=h).\tam Тогда

\bdf
T^{\ad(a)}(\xi)=T^{\ad(a)}\left(\frac{\alpha\mapr{}\beta}{g}\right)=
A^{g\alpha=\beta g}_{g}=\delta^{ag}_{g\alpha=\beta g}-\delta^{ga}_{g\alpha=\beta g}.
\edf

Первое слагаемое функции \mat T^{\ad(a)}(\xi)\tam  равно единице тогда и только тогда, когда \blf \beta=a\elf, т.е. когда морфизм
\blf\xi\in \Mor(g^{-1}ag,a)\elf. Аналогично, второе слагаемое
функции \mat T^{\ad(a)}(\xi)\tam  равно минус единице  тогда и только тогда, когда \blf \alpha=a\elf,
т.е. когда морфизм
\blf\xi\in \Mor(a,gag^{-1})\elf.

Другими словами, на морфизмах \blf \Mor(g^{-1}ag,a)\elf матрица
\blf\ad(a)\elf равна единице,  на морфизмах \elf \Mor(a,gag^{-1})\elf
равна минус единице и на морфизмах \elf \Mor(u,u)\elf равна нулю, равно как и на морфизмах \elf \Mor(a,a)\elf и \elf \Mor(v,v)\elf в следующей диаграмме:

\bdf
\xymatrix{
u\ar@(lu,ru)[]^{T=0}\ar[r]\ar@/_/[r]_{T=+1}
&a\ar@/_/[r]_{T=-1}\ar[r]
\ar@(lu,ru)[]^{T=0}&v\ar@(lu,ru)[]^{T=0}
}
\edf

Отсюда получаем теорему:
\begin{theoremru}[О внутренних деривациях]
\end{theoremru}
Характеры внутренних дериваций тривиальны на \mat\Mor(a,a)\tam:
\bdf
\xymatrix{
\Int(\cA)\ar[r]\ar[d]_{\cap}&\Der(\cA)\ar[d]^{T}_{\approx}\\
\ker p_{a}\ar[r]& \T_{f}(\cG)\ar[r]^{p_{a}}&
\T_{f}(\Mor(a,a))
}
\edf
Множество \mat\T_{f}(\Mor(a,a))\tam совпадает с группой всех характеров
\bdf
\T_{f}(\Mor(a,a))=\T(\Mor(a,a)).
\edf

Заметим, что если характер \mat T\in \T(\cG)\tam  равен нулю на \mat\Mor(a,a),\tam
то он равен нулю на всяком \mat\Mor(u,u)\tam для сопряженного элемента \mat u\in[a], \quad u=gag^{-1}.\mat

Так что диаграмма имеет следующий вид:
\bdf
\xymatrix{
\Int(\cA)\ar[r]\ar[d]_{\cap}^{T}&\Der(\cA)\ar[d]^{T}_{\approx}\\
\ker p_{a}\ar[r]& \T_{f}(\cG)\ar[r]^{p_{a}}&
\T(\Mor(a,a))
}
\edf

С точки зрения деривационной проблемы Джонсона (Johnson(2001)\cite{Johnson-2001}) о деривациях естественно обозначить через
\mat\Out(\cA)\tam факторгруппу \mat\Out(\cA)=\Der(\cA)/\Int(\cA)\tam и назвать ее алгеброй внешних дериваций алгебры \mat\cA\tam. Таким образом предыдущая диаграмма дополняется до диаграммы
\bdf
\xymatrix{
0\ar[r]&\Int(\cA)\ar[r]\ar[d]_{\cap}^{T}&\Der(\cA)\ar[d]^{T}_{\approx}\ar[r]&\Out(\cA)\ar[d]\ar[r]&0\\
0\ar[r]&\ker p_{a}\ar[r]& \T_{f}(\cG)\ar[r]^{p_{a}}&
\T(\Mor(a,a))
}
\edf

\subsection*{Описание внутренних дериваций:}

Прежде всего заметим, что множество морфизмов \mat\Mor(\cG)\tam группоида \mat\cG\tam распадается на несвязное объединение морфизмов по классам сопряженных элементов группы \mat G,\tam являющихся по определению объектами
группоида \mat\cG.\tam Сама группа \mat G\tam разлагается в несвязное объединение классов сопряженных элементов
\bdf
G=\coprod\limits_{g\in G}[g], \quad [g]=\{h: \exists a\in G, h=aga^{-1}\}.
\edf
Соответственно, множество морфизмов тоже представляется в виде несвязного объединения
\bdf
\Mor(\cG)=\coprod\limits_{[g]}\Mor(\cG_{[g]}).
\edf
Это значит, что построение любой деривации можно строить в виде дериваций \mat\Der_{[g]}(\A)\tam независимо в каждой подкатегории \mat\cG_{[g]}\mat \ как локально финитные характеры на каждой из них.

Естественная проблема заключается в том, чтобы установить, верно ли, что
множество всех дифференцирований,  тривиальных  на всех \blf\Mor(u,u),\elf совпадает с множеством внутренних дифференцирований. Другими словами, верно ли, что вложение
\mat\xymatrix{
\Int_{\hskip -4pt[u]}(\cA)\ar[r]^{\subset}&
\ker p_{u}}\tam является изоморфизмом:
\bdf\xymatrix{
\Int_{\hskip -4pt [u]}(\cA)\ar[r]^{\subset}&
\ker p_{u}\ar[r]^{\subset} &\T_{f}(\cG)\ar[r]^{p_{u}}&
\T(\Mor(u,u))
} \edf

Анализ этой проблемы позволяет сформулировать конкретные условия на локально финитный характер \mat T:\Mor(\cG)\mapr{}R\tam, который реализует внутреннюю деривацию \mat X\in\Int(\cG)\tam,
\bdf
T=T^{X}.
\edf

\subsubsection*{Случай единичного элемента \mat[e]:\tam}

В частности, одна из подкатегорий соответствует единичному элементу \mat e\in G,\tam у которого \mat[e]=\{e\}\tam. В этом частном случае подкатегория \mat\cG_{[e]}\tam состоит из одного объекта \mat e\in G\tam, а множество морфизмов изоморфно группе \mat G,\tam \mat\Mor(e,e)\approx G.\tam В частности, множество локально финитных характеров \mat T_{f}(\cG_{[e]})\approx T(\cG_{[e]})\approx T(G)\tam
изоморфно группе всех характеров на группе \mat G.\tam Каждый характер на группе \mat T\in T(G)\tam реализуется деривацией \mat X\in\Der(G),\tam \mat T^{X}=T\tam. Действительно, характер \mat T\in T(G)\tam -- это характер на категории \mat\cG\tam, который равен \mat T\tam на \mat\Mor(e,e)\tam и равен нулю на всех остальных слагаемых \mat\Mor(\cG_{[g]}), \quad g\neq e\tam. Поэтому соответствующая матрица \mat \| X^{h}_{g}\|\tam оператора \mat X\tam задается формулой
\bdf
X^{h}_{g}=T(g)\delta^{h}_{g}.
\edf
Все деривации, соответствующие характерам на подкатегории \mat\cG_{[e]}\tam, не являются внутренними деривациями.

Аналогичные рассуждения годятся и для других классов сопряженности, состоящих из конечного числа элементов, т.е. когда \mat\#[g]<+\infty.\mat \ В частности,
это справедливо для элементов из центра \mat g\in Z(G).\tam

\subsubsection*{Точная последовательность:}
Если отказаться от ограничения на характеры локальной финитности, то можно установить точность последовательности, как это сформулировано в следующей теореме:

\begin{theoremru}\label{t3}
Последовательность
\bdf
0\mapr{}\ker p_{a}\mapr{}\T(\cG_{[a]})\mapr{p_{a}}\T(\Mor(a,a))\mapr{}0, \quad a\in U,
\edf
точна.
\end{theoremru}

\proofru
Надо проверить только эпиморфность отображения \blf p(a)\elf.

Пусть \blf\chi\in\X(G^{a})\elf.

Фиксируем элемент \blf a\in U\elf и фиксируем произвольные элементы \blf g_{a,b}\in\Mor(a,b),\forall b\in U, \elf удовлетворяющие
условию \blf g_{a,a}=e\elf. Такие элементы существуют, поскольку \blf\Mor(a,b)\elf не пусто. Далее полагаем
\bdf
g_{b,c}=g_{a,c}g_{a,b}^{-1}\in\Mor(b,c).
\edf
Заметим, что условие \blf g_{b,c}\in\Mor(b,c)\elf означает, что
\bdf
g_{b,c}b=c g_{b,c}.
\edf
Проверим это условие. Во-первых,  по построению
\bdf
g_{a,b}a=b g_{a,b}.
\edf
Далее, \blf g_{b,a}=g_{a,a}g_{a,b}^{-1}=g_{a,b}^{-1}\in\Mor(b,a),\elf
т.е.
\bdf
g_{b,c}=g_{a,c}g_{b,a}.
\edf
Значит,
\bdf
g_{b,c}b=g_{a,c}g_{b,a}b=g_{a,c}ag_{b,a}=
c g_{a,c}g_{b,a}=c g_{b,c},
\edf
т.е.
\bdf
g_{b,c}\in\Mor(b,c).
\edf

Следовательно, имеет место соотношение
\bdf
g_{c,d}g_{b,c}=g_{b,d}\in\Mor(b,d).
\edf

В самом деле
\bdf
g_{c,d}g_{b,c}=g_{a,d}g_{c,a}g_{a,c}g_{b,a}=
g_{a,d}g_{b,a}=g_{b,d}\in\Mor(b,d).
\edf

Теперь построим характер \blf X\elf  на группоиде \blf\cG_{U},\elf
\blf X:\Mor(U)\mapr{}\R\elf.
Пусть \blf x_{b,c}\in\Mor(b,c)\elf -- произвольный морфизм. Тогда
\bdf
g_{c,a}x_{b,c}g_{a,b}\in\Mor(a,a)=G^{a}.
\edf
Полагаем
\bdf
X(x_{b,c})\bydef\chi(g_{c,a}x_{b,c}g_{a,b}).
\edf
Легко проверить, что отображение \blf X:\Mor(U)\mapr{}\R\elf аддитивно:
\bdf
\begin{array}{l}
X(x_{c,d}x_{b,c})=\chi(g_{d,a}(x_{c,d}x_{b,c})g_{a,b})=
\chi(g_{d,a}x_{c,d}g_{a,c}g_{c,a}x_{b,c}g_{a,b})=\\\\=
\chi(g_{d,a}x_{c,d}g_{a,c})\chi(g_{c,a}x_{b,c}g_{a,b})=
X(x_{c,d})X(x_{b,c}).
\end{array}
\edf

Ограничение характера \blf X\elf на \blf\Mor(a,a)=G^{a}\elf совпадает с \blf\chi\elf:
\bdf
X(x_{a,a})=\chi(g_{a,a}x_{a,a}g_{a,a})=\chi(x_{a,a}).
\edf
\qed

\subsubsection*{Редукция к группам коцепей}
Вернемся к изучению
характеров \mat T,\tam тривиальных на подпространстве \mat\Mor(a,a)\tam
\bdf\xymatrix{
\Int_{a}\cA\ar[r]^{\subset}&
\ker p^{f}_{a}\ar[r]^{\subset} &\T_{f}(\cG_{[a]})\ar[r]^{p^{f}_{a}}&
\T(\Mor(a,a))
} \edf
Обозначим через \mat\Delta(\cG_{[a]})\tam симплекс, вершинами которого являются элементы класса сопряженности \mat[a]\subset G=\Obj(\cG_{[a]}).\tam

Поскольку каждый характер \mat T\in\ker p^{f}_{a}\tam принимает одинаковое значения на множестве всех морфизмов
\mat\Mor(b,c)\tam, \mat b,c\in[a]\subset G\tam, то это значит, что существует естественное вложение
\bdf
\varphi:\ker p^{f}_{a}\hookrightarrow \C^{1}(\Delta(\cG_{[a]})),
\edf
в группу коцепей симплекса \mat\Delta(\cG_{[a]})\tam, причем каждый характер \mat T\in\ker p^{f}_{a}\tam
является коциклом:
\bdf
\xymatrix{
\C^{0}_{f}(\Delta(\cG_{[a]}))\ar@{^{(}->}[r]&\C^{0}(\Delta(\cG_{[a]}))\ar[r]^{\delta}&\C^{1}(\Delta(\cG_{[a]}))\ar[r]^{\delta}&C^{2}(\Delta(\cG_{[a]}))\\
\Int_{a}(\cA)\ar@{^{(}->}[r]^{T}\ar@{^{(}->}[ur]\ar@{^{(}->}[u]^{\varphi_{i0}}&
\ker p^{f}_{a}\ar@{^{(}->}[u]^{\varphi_{0}}\ar@{^{(}->}[ur]^{\varphi}\ar[r]&0\ar[ur]
}
\edf

Вложение \mat \varphi_{i0}\tam является изоморфизмом. Образ
\bdf
\delta(\varphi_{i0}(\Int_{a}(\cA)))\subset
\delta(C^{0}_{f}(\Delta(\cG_{[a]})))\subset
\ker p^{f}_{a}\subset
\C^{1}(\Delta(\cG_{[a]}))
\edf
может быть описан как некоторое множество коциклов, удовлетворяющих определенным условиям. Рассмотрим класс сопряженности \mat[a],\tam на котором действует группа \mat G\tam при помощи присоединенного действия
\bdf
G\times[a]\mapr{\ad}[a], \quad \ad_{g}(b)=gbg^{-1}, \, b\in[a].
\edf
Для произвольного элемента \mat g\in G\tam рассмотрим граф \mat\Gamma_{g,a}\subset \Delta(\cG_{[a]}),\tam составленный из ориентированных ребер, начало которых являются элементы
\mat b\in[a],\tam а концы равны
\mat gbg^{-1}\in[a]\tam. Таким образом, каждое ребро имеет вид направленного отрезка
\mat[b, gbg^{-1}].\tam Граф \mat\Gamma_{g,a}\tam разлагается в несвязное объединение направленных путей, составленных из направленных ребер
\bdf
\Gamma_{g,a}=\coprod\limits_{\alpha}\Gamma_{g,a}^{\alpha}.
\edf
Каждый направленный путь \mat\Gamma_{g,a}^{\alpha}\tam может быть бесконечным в обе стороны или конечным, и тогда является циклическим путем.

\begin{theoremru}
Пусть \mat X\in\Int_{[a]}(\cA).\tam Тогда коцепь
\mat \varphi(T^{X})\in\C^{1}(\Delta(\cG_{[a]}))\tam удовлетворяет условию:
\begin{itemize}
\item[\F\F:]
Коцепь \mat \varphi(T^{X})\tam
является финитной на графе \mat\Gamma_{g,a},\tam а сумма значений коцепи
\mat \varphi(T^{X})\tam на каждом направленном пути \mat\Gamma_{g,a}^{\alpha}\tam
равна нулю.
\end{itemize}\end{theoremru}

Неясно, верно ли обратное утверждение, т.е если оператор \mat X\in\Der(\cA),\tam \mat T^{X}\in\ker p^{f}_{a}\subset T_{f}(\cG_{[a]})\tam, удовлетворяет условию (\F\F), верно ли что
\mat
X\in\Int_{[a]}(\cA)?
\tam

\subsubsection*{Ядро \mat\ker p^{f}_{a}\tam не совпадает с внутренними деривациями}

Существуют примеры групп, для которых имеются локально финитные характеры
\mat T\in \ker p^{f}_{a}\tam, которые не удовлетворяют условию (\F\F).
Рассмотрим  в качестве простейшего примера свободную группу с двумя образующими
\mat G=F<x_{1},x_{2}>.\tam В качестве класса сопряженных элементов возьмем класс \mat[x_{1}]\subset G.\tam Рассмотрим характер \mat T:\Mor(\cG_{[x_{1}]})\mapr{} R\tam, задавая значения характера на образующих группоида \mat\Mor(\cG_{[x_{1}]})\tam независимо друг от друга. Множество образующих группоида \mat\Mor(\cG_{[x_{1}]})\tam состоит из морфизмов вида
\bdf
\xi=\left(\frac{\alpha\mapr{}\beta}{g}\right), \quad g=x_{1},x_{2}, \quad \beta=g\alpha g^{-1}, \quad \alpha\in[x_{1}].
\edf
Полагаем
\bdf
\begin{array}{l}
T\left(\frac{x_{2}x_{1}x_{2}^{-1}\mapr{}x_{1}x_{2}x_{1}x_{2}^{-1}x_{1}^{1}}{x_{1}}\right)=1;\\
\hbox{на остальных образующих характер }T\hbox{ равен нулю}.
\end{array}
\edf
Поскольку в свободной группе нет соотношений кроме естественных сокращений в словах, то функцию \mat T\tam можно по аддитивности продолжить до некоторого характера на группоиде \mat\cG_{[x_{1}]}.\tam

Следовательно, на одном из направленных путей вида
\mat\Gamma^{\alpha}_{x_{1},x_{1}}\tam характер \mat Е\tam имеет только одно значение, равное единице, а на всех остальных направленных путях этот характер тождественно равен нулю. Это означает, что условие (\T1) выполнено, но условие (\F\F) не выполняется.

\section{Дополнение: группы с конечным числом образующих}

Пусть группа $G$ конечно представима, т.е. имеет конечное число образующих $\{g_{1}, g_{2}, \dots, g_{k}\}$ и конечное число определяющих
соотношений $\{S_{1}, S_{2}, \dots, S_{l}\}$:
$$
G=F<g_{1}, g_{2}, \dots, g_{k}>/\{S_{1}, S_{2}, \dots, S_{l}\}.
$$

Каждое соотношение \mat S_{i}\tam -- это слово длины \mat s_{i}\tam, составленное из образующих элементов \mat g_{1}, g_{2}, \dots, g_{k}\tam или обратных к образующим элементов \mat g_{1}^{-1}, g_{2}^{-1}, \dots, g_{k}^{-1}\tam. Так что соотношение \mat S_{i}\tam представляется в следующем виде
\bdf
S_{i}=g_{j_{(i,1)}}^{\varepsilon_{(i,1)}}
g_{j_{(i,2)}}^{\varepsilon_{(i,2)}}
g_{j_{(i,3)}}^{\varepsilon_{(i,3)}}
\dots
g_{j_{(i,s_{i})}}^{\varepsilon_{(i,s_{i})}}\equiv 1,
\quad \varepsilon_{(i,j)}=\pm 1.
\edf

Каждое соотношение \mat S_{i}\tam индуцирует целую серию соотношений на группоиде \mat\cG\tam:
\bdf
\xi_{(i,1)}\xi_{(i,1)}\xi_{(i,3)}\dots\xi_{(i,s_{i})}\equiv
\left(\frac{\alpha_{(i,1)}\mapr{}\alpha_{(i,1)}}
{1}\right),
\edf
где морфизмы \mat \xi_{(i,j)}\tam определяются по правилу
\bdf
\begin{array}{l}
\xi_{(i,j)}=\left(\frac{\alpha_{(i,j)}\mapr{}\beta_{(i,j)}}
{g_{j_{(i,j)}}^{\varepsilon_{(i,j)}}}\right),\\\\
\beta_{(i,j)}=g_{j_{(i,j)}}^{\varepsilon_{(i,j)}}\alpha_{(i,j)}
g_{j_{(i,j)}}^{-\varepsilon_{(i,j)}},\\\\
\beta_{(i,j)}=\alpha_{(i,j+1)},\quad
\beta_{(i,s_{i})}=\alpha_{(i,1)}.
\end{array}
\edf

Таким образом, чтобы задать локально финитный характер \mat T\tam на группоиде \mat\cG\tam, достаточно определить значения
локально финитного характера \mat T\tam на множестве образующих
\bdf
\coprod\limits_{i=1}^{k}\cH_{g_{i}},
\edf
таким образом, чтобы на каждом соотношении
\bdf
\xi_{(i,1)}\xi_{(i,1)}\xi_{(i,3)}\dots\xi_{(i,s_{i})}\equiv
\left(\frac{\alpha_{(i,1)}\mapr{}\alpha_{(i,1)}}
{1}\right),
\edf
выполнялось условие аддитивности:
\bdf
T(\xi_{(i,1)})+T(\xi_{(i,1)})+T(\xi_{(i,3)}+\cdots+T(\xi_{(i,s_{i})})=0
\edf


\end{otherlanguage}

\end{document}